\documentclass[12pt,reqno,oneside]{amsart}
\usepackage[T1]{fontenc}
\usepackage{lmodern}
\usepackage{amssymb}
\usepackage[letterpaper,textwidth=5.95in,textheight=8.7in,centering]{geometry}
\usepackage[expansion=false]{microtype}
\usepackage{needspace}
\usepackage[hidelinks]{hyperref}

\hypersetup{
  pdftitle={A mean index criterion for multiple closed characteristics},
  pdfsubject={A mean index criterion for closed characteristics in all half-dimensions},
  pdfkeywords={closed characteristic, mean index, common index jump, equivariant Morse theory}
}
\newtheorem{theorem}{Theorem}[section]
\newtheorem{lemma}[theorem]{Lemma}
\newtheorem{corollary}[theorem]{Corollary}
\theoremstyle{remark}
\newtheorem{remark}[theorem]{Remark}
\newtheorem*{unnumremark}{Remark}
\newtheorem{example}[theorem]{Example}
\numberwithin{equation}{section}
\newcommand{\R}{\mathbb{R}}
\newcommand{\Z}{\mathbb{Z}}
\newcommand{\N}{\mathbb{N}}
\newcommand{\PP}{\mathcal{P}}
\newcommand{\mchi}{\widehat{\chi}}

\title[A mean index criterion for closed characteristics]
{A mean index criterion for multiple closed characteristics}

\author{Zhiping Fan}
\address{School of Mathematical Sciences, Zhejiang Normal University,
Jinhua, Zhejiang 321004, China}
\email{fanzhiping@mail.nankai.edu.cn}
\thanks{Zhiping Fan was partially supported by NSFC (No.~12201566)
and Zhejiang Normal University (YS304022950).}
\author{Zhihao Zhao\textsuperscript{*}}
\address{School of Mathematics, China University of Mining and Technology,
Xuzhou, Jiangsu 221000, China}
\email{zhaozh@cumt.edu.cn}
\thanks{Zhihao Zhao was partly supported by the NSFC (No.~12401129).}
\thanks{\textsuperscript{*}Corresponding author.}
\date{}

\begin{document}

\begin{abstract}
Let \(n\ge2\). We prove that a compact strictly star-shaped \(C^3\) hypersurface in \(\mathbb R^{2n}\) carries at least \(n\) geometrically distinct prime closed characteristics if all iterates are nondegenerate and all mean indices are at least one. It suffices that all mean indices are positive and that the bound holds only for prime characteristics with negative mean Euler characteristic in the variational grading. The proof uses a minimum argument for the global Morse remainder and requires no individual index exclusion. For odd \(n\ge3\), we obtain a smaller, dimension-dependent strict bound when good iterates of Maslov-type index \(-1\) are excluded. Irrational ellipsoids show that the multiplicity bound is sharp.
\end{abstract}

\maketitle

\section{Introduction}\label{sec:intro}

Let \(\Sigma\subset\R^{2n}\) be a compact \(C^3\) hypersurface
strictly star-shaped with respect to the origin. Thus each ray
from the origin meets \(\Sigma\) once and transversely. The
standard Liouville form
\(\lambda_0=\tfrac12\sum_{j=1}^n(x_j\,dy_j-y_j\,dx_j)\)
restricts to a contact form \(\lambda\) on \(\Sigma\).
A closed characteristic of period \(\tau>0\) is a solution of
\[
 \dot y=R_\lambda(y),\qquad y(t+\tau)=y(t),
\]
where \(\lambda(R_\lambda)=1\) and
\(d\lambda(R_\lambda,\cdot)=0\). It is \emph{prime} if \(\tau\)
is its least period. We identify solutions that differ by a time
translation; two prime characteristics are geometrically distinct
if their images differ. The multiplicity problem asks whether every
such hypersurface has at least \(n\) geometrically distinct prime
closed characteristics. We give a sufficient mean index condition
for this conclusion in every half-dimension \(n\ge2\).

Write \(i(y,m)\) for the Maslov-type index of the \(m\)-fold
iterate \(y^m\), computed in the ambient symplectic trivialization
using \(j_\Sigma^\alpha\), \(1<\alpha<2\), where
\(j_\Sigma(rx)=r\) for \(x\in\Sigma\) and \(r>0\).
The transverse Conley--Zehnder index, computed in a trivialization
extending over a capping disk, is
\begin{equation}\label{eq:grading}
 \mu_y(m)=\mu_{\mathrm{CZ}}(y^m)=i(y,m)+1.
\end{equation}
The conventions and this shift are explained in
Section~\ref{sec:prelim}. The mean index is
\[
 \Delta_y=\lim_{m\to\infty}\frac{i(y,m)}m
         =\lim_{m\to\infty}\frac{\mu_y(m)}m.
\]
An iterate is \emph{nondegenerate} if its transverse return map
has no eigenvalue one. Write \(\N=\{1,2,\ldots\}\).
For a prime characteristic whose iterates
are all nondegenerate, set
\[
 G_y=\{m\in\N:i(y,m)\equiv i(y,1)\pmod2\}.
\]
We call \(y^m\) \emph{good} if \(m\in G_y\), and \emph{bad}
otherwise. The mean Euler characteristic of the prime orbit \(y\),
in the variational grading, is
\begin{equation}\label{sec2:2}
 \mchi(y)=\lim_{r\to\infty}\frac1r
      \sum_{\substack{1\le m\le r\\m\in G_y}}
      (-1)^{\mu_y(m)-(n+1)}.
\end{equation}
The limit exists and belongs to \(\{1,1/2,-1,-1/2\}\), as recalled
in Section~\ref{subsec:morse}. The exponent is the variational
Morse degree \(I(y^m)=\mu_y(m)-(n+1)\). Thus \eqref{sec2:2}
uses the convention of~\cite[Definition~4.8 and Remark~4.9]{LLW}
for both parities of \(n\).

\begin{theorem}\label{thm:main}
Let \(n\ge2\), and let \(\Sigma\subset\R^{2n}\) be a
compact \(C^3\) hypersurface strictly star-shaped with respect to
the origin. Assume that
\begin{itemize}
 \item[\textup{(H1)}] every iterate of every prime closed
 characteristic is nondegenerate;
 \item[\textup{(H2)}] every prime closed characteristic \(y\)
 has \(\Delta_y>0\); if \(\mchi(y)<0\), then \(\Delta_y\ge1\).
\end{itemize}
Then \(\Sigma\) carries at least \(n\) geometrically distinct
prime closed characteristics.
\end{theorem}

\begin{corollary}\label{cor:uniform}
Let \(n\ge2\), and let \(\Sigma\subset\R^{2n}\) be a
compact strictly star-shaped \(C^3\) hypersurface. If every
iterate is nondegenerate and every prime closed characteristic
has mean index at least one, then \(\Sigma\) carries at least
\(n\) geometrically distinct prime closed characteristics.
\end{corollary}

There is no restriction on individual Maslov-type or
Conley--Zehnder indices in these statements. In particular, good
iterates of Conley--Zehnder index zero are permitted. The numerical
bound in \textup{(H2)} concerns only the characteristics whose good
iterates have Conley--Zehnder index congruent to \(n\) modulo two.
Irrational ellipsoids have
exactly \(n\) prime characteristics and satisfy the hypotheses;
see Example~\ref{ex:ellipsoid}. This proves sharpness of the orbit count. It remains unclear whether one is the least possible mean index
bound.

There is also a rigidity statement when the lower bound $n$ is attained.
Under \textup{(H1)}, assume that all mean indices are positive and
that \(\Delta_y>1\) whenever \(\mchi(y)<0\).
If there are exactly \(n\) prime closed characteristics, all iterates
are good and the Reeb flow is dynamically convex; see
Corollary~\ref{cor:minimal-rigidity}.

The coefficient estimate also yields a smaller bound if an index
exclusion is imposed. For odd \(n\ge3\), put
\begin{equation}\label{eq:secondary-threshold}
 D_{\max}(n)=\left\lfloor\frac{2n-3}{3}\right\rfloor,
 \qquad \theta_n=1-\frac1{D_{\max}(n)}.
\end{equation}

\Needspace{11\baselineskip}
\begin{theorem}\label{thm:smaller}
Let \(n\ge3\) be odd, and let \(\Sigma\subset\R^{2n}\) be a
compact strictly star-shaped \(C^3\) hypersurface satisfying
\textup{(H1)}. Suppose that every prime characteristic has positive
mean index and that
\begin{itemize}
 \item[\textup{(H3)}] every good iterate \(y^m\) satisfies
 \(i(y,m)\ne-1\), equivalently \(\mu_y(m)\ne0\).
\end{itemize}
If \(\Delta_y>\theta_n\) whenever \(\mchi(y)<0\), then
\(\Sigma\) has at least \(n\) geometrically distinct prime closed
characteristics.
\end{theorem}

For \(n=3,5,7,9,11\), the values of \(\theta_n\) are,
respectively, \(0,1/2,2/3,4/5,5/6\). The inequality in
Theorem~\ref{thm:smaller} is strict and its endpoint is not asserted.
For \(n=3\), positivity alone supplies the mean index condition.
This case is already covered by Duan, Liu, Long, Qi, and
Wang~\cite{DLLQW}: they require only that prime characteristics
avoid Maslov-type index \(-1\), without any mean index assumption.

For even \(n\), \textup{(H1)}, \textup{(H3)}, and nonvanishing of
all mean indices already imply the lower bound \(n\), by
\cite[Definition~1.1 and Theorem~1.2]{DLLW}. No quantitative lower
bound on the mean indices is needed in that case. We therefore
state Theorem~\ref{thm:smaller} only for odd \(n\).

For smooth hypersurfaces we can also allow negative mean indices.
This follows from Theorem~\ref{thm:main} and the negative-mean-index
criterion of \c{C}ineli, Ginzburg, and G\"urel
\cite[Theorem~B]{CGGFiltered}.

\begin{corollary}\label{cor:absolute}
Let \(n\ge2\), and let \(\Sigma\subset\R^{2n}\) be a
compact strictly star-shaped \(C^\infty\) hypersurface. If every
iterate is nondegenerate and \(|\Delta_y|\ge1\) for every prime
closed characteristic, then \(\Sigma\) carries at least \(n\)
geometrically distinct prime closed characteristics.
\end{corollary}

The existence theory starts with Rabinowitz's theorem on periodic
solutions on star-shaped energy surfaces~\cite{Rab}. Berestycki,
Lasry, Mancini, and Ruf~\cite{BLMR} obtained multiplicity under
pinching between homothetic ellipsoids. Viterbo~\cite{Viterbo}
developed equivariant Morse theory for star-shaped Hamiltonian
systems. Long's iteration formula~\cite{Long} and the common index
jump theorem of Long and Zhu~\cite{LZ} made it possible to compare
the contributions of different prime orbits to these Morse counts.
For strictly convex hypersurfaces, Long and Zhu obtained
\(\lfloor n/2\rfloor+1\) closed characteristics in general and
\(n\) in the nondegenerate case. For nondegenerate star-shaped
hypersurfaces, Hu and Long~\cite{HL} proved the existence of at
least two, without a mean index assumption. Thus the multiplicity
conclusion of Theorem~\ref{thm:main} is already known when \(n=2\).
We include this case because the proof applies to it without a
separate argument. The full bound \(n\) under dynamical convexity, including
degenerate flows, was established by \c{C}ineli, Ginzburg, and
G\"urel~\cite[Theorem~A]{CGG}.

Another class of results excludes selected low-index iterates.
For odd \(n\), Duan, Liu, Long, and Wang
\cite[Theorem~1.2]{DLLW} obtained \(n\) closed characteristics
under nonzero mean indices and exclusion of the good Maslov-type
indices \(-2,-1,0\). Duan and Xie~\cite[Theorem~1.1]{DX}
obtained \((n+3)/2\) under nonzero mean indices and exclusion of
the Maslov-type index \(-1\) on all iterates. For \(n=5\),
Duan and Fan~\cite[Theorem~1.1(a)]{DF} proved the bound five
when all mean indices are at least one and good iterates avoid
\(-1\). Corollary~\ref{cor:uniform} removes that index exclusion
and applies in every half-dimension \(n\ge2\).

The corresponding contact formulation uses \eqref{eq:grading}.
The radial parametrization \(F(u)=f(u)u\), \(u\in S^{2n-1}\),
satisfies \(F^*\lambda=f^2\lambda_0|_{S^{2n-1}}\), so the induced
contact structure is the standard one. Under nonzero mean indices
and exclusion of good Conley--Zehnder index zero, Liu
\cite[Corollary~1.8]{Liu26} and Xie
\cite[Theorem~1.5(ii)]{Xie26} give \((n+3)/2\) prime orbits on
\(S^{2n-1}\) for odd \(n\). In Xie's formula, the base
\(\mathbb{CP}^{n-1}\) has Euler characteristic \(n\) and middle
Betti number one; the lower bound in~\cite[(1.4)]{Xie26}, at
\(p=0\), is
\[
 \max\left\{\frac{n+3}{2},\,2-n\right\}=\frac{n+3}{2}.
\]
Abreu and Macarini~\cite{AM} obtain multiplicity bounds from
periodic equivariant symplectic homology; their sharp bound on
these spheres uses other index exclusions. Theorem~\ref{thm:main}
assumes a quantitative positive mean index bound on the prime
characteristics with negative mean Euler characteristic and imposes
no individual index exclusion. Thus it
does not replace the conclusions of those results for arbitrary
nonzero mean indices.

The proof uses the precise iteration formula and the enhanced
common index jump theorem~\cite{DLW,DLLW}. Counts in reflected
degrees and complementary choices of jumps appear in Duan and
Fan~\cite[\S4]{DF}. Alternating sums related by index recurrence
are also used by Abreu and Macarini~\cite[\S6]{AM}; see
Xie~\cite[Remark~1.7]{Xie26}. The two boundary terms of an
alternating Morse sum are already present in
Viterbo~\cite[(8.6)--(8.7)]{Viterbo}.
The additional step in our proof is a minimum argument for the
global Morse remainder. Assuming that there are finitely
many prime closed characteristics, write the Morse relation as
\(C_j=B_j+U_j+U_{j-1}\), where \(C_j\) counts good iterates
of degree \(j\), \(B_j\) is the corresponding Betti number,
and \(U_j\) is a fixed sequence of nonnegative integers,
independent of the common index jump.
For every integer \(p\ge n\) with \(p\equiv n\pmod2\), a
sufficiently large even common
index jump center \(d\) gives, after retaining the boundary terms
in the alternating Morse sums, the identity
\[
 U_{d+p-1}+U_{d-p}=2U_{p-1}-D,
 \qquad D=n-a+c,
\]
where \(a\) and \(c\) count the prime characteristics with mean
Euler characteristic \(1\) and \(-1\), respectively.
Nonnegativity alone gives \(U_{p-1}\ge D/2\). To obtain the bound
needed here, first choose \(p_0\ge n\), with
\(p_0\equiv n\pmod2\), at which \(U_{p-1}\) attains its minimum
over this range. Then choose the common index jump for this
\(p_0\). Since \(U_{d+p_0-1}\ge U_{p_0-1}\) and
\(U_{d-p_0}\ge0\), the identity gives \(U_{p_0-1}\ge D\).
Consequently \(C_p\ge U_{p-1}\ge D\) throughout the range.
These bounds give a lower bound on the asymptotic density of good
iterates in degrees congruent to \(n\) modulo two. Counting these
iterates separately on their prime orbits expresses the same density
in terms of the mean indices. Under \textup{(H2)}, the two bounds
are incompatible with fewer than \(n\) prime characteristics.
Thus, after choosing the minimizing degree, one application of the
common index jump theorem suffices for the main argument.
Complementary jumps and the resonance identity are used in the proof
of Theorem~\ref{thm:smaller}.

Section~\ref{sec:prelim} fixes the index conventions and the
global Morse coefficients. Section~\ref{sec:jump} constructs a
common index jump. Section~\ref{sec:counting} gives the minimum argument and proves
the asymptotic counting formula. It also derives a quantitative
multiplicity bound and a rigidity statement when the minimal
orbit count is attained.
Section~\ref{sec:main-proof} proves the main result and its
corollaries. Section~\ref{sec:smaller} proves
Theorem~\ref{thm:smaller} and discusses its endpoint.

\noindent\emph{Notation.}
For \(x\in\R\),
\(\lfloor x\rfloor\) and \(\lceil x\rceil\) denote the floor
and ceiling, and \(\{x\}=x-\lfloor x\rfloor\). Set
\[
 \operatorname{rnd}(x)=\left\lfloor x+\tfrac12\right\rfloor,
 \qquad \|x\|_{\Z}=|x-\operatorname{rnd}(x)|
                  =\min_{h\in\Z}|x-h|.
\]
Thus \(\operatorname{rnd}\) assigns half-integers to the larger
integer.

\section{Preliminaries}\label{sec:prelim}

Throughout Sections~\ref{sec:prelim}--\ref{sec:counting},
\(n\ge2\), \(\Sigma\) is as in Theorem~\ref{thm:main},
and \textup{(H1)} holds. We assume that the set \(\PP\) of
geometrically distinct prime closed characteristics is finite
and that \(\Delta_y>0\) for every \(y\in\PP\).
The existence theorem~\cite{Rab} ensures that \(\PP\ne\varnothing\).
Write
\[
 q=\#\PP,\qquad \delta_*=\min_{y\in\PP}\Delta_y>0.
\]
The multiplicity conclusions are immediate when \(\PP\) is
infinite, so it suffices to consider the finite case. No numerical lower
bound on \(\Delta_y\), and no index exclusion, is imposed in
these three sections.

\Needspace{9\baselineskip}
\subsection{The linearized flow and its indices}\label{subsec:indices}

Put \(\omega_0=d\lambda_0\). The radial field \(Y(z)=z/2\) satisfies
\(\iota_Y\omega_0=\lambda_0\) and is transverse to \(\Sigma\).
Consequently \(\lambda\) is a contact form and \(d\lambda\) is
symplectic on \(\xi=\ker\lambda\). The Reeb flow \(\phi_\lambda^t\)
preserves \(\lambda\). For a prime characteristic \(y\) of least Reeb
period \(\tau_y\), let
\[
 P_y=d\phi_\lambda^{\tau_y}(y(0))|_{\xi_{y(0)}}.
\]
The iterate \(y^m\) traverses the same orbit \(m\) times in time
\(m\tau_y\) and has transverse return map \(P_y^m\). Thus
\textup{(H1)} is equivalent to the absence of roots of unity in the
spectrum of \(P_y\) for every \(y\in\PP\).

In the coordinates \(z=\sqrt r\,x\), where \(r>0\) and \(x\in\Sigma\),
we have \(\lambda_0=r\lambda\). The Hamiltonian \(H=j_\Sigma^\alpha\)
is \(h(r)=r^{\alpha/2}\). With \(\iota_{X_H}\omega_0=-dH\), its
flow is
\[
 \phi_H^t(r,x)=\bigl(r,\phi_\lambda^{\,t h'(r)}(x)\bigr).
\]
One traversal of \(y\) on \(r=1\) therefore takes Hamiltonian time
\(T_y=\tau_y/h'(1)=2\tau_y/\alpha\). In the standard trivialization of
\(T\R^{2n}\) put \(\gamma_y(t)=d\phi_H^t(y(0))\) for
\(0\le t\le T_y\), and set \(M_y=\gamma_y(T_y)\).
The ambient linearized path for \(y^m\) is the path of \(m\) successive
periods
\[
 \gamma_y^{[m]}(t)=\gamma_y(t-jT_y)M_y^j,\qquad
 jT_y\le t\le(j+1)T_y,\quad j=0,\ldots,m-1 .
\]
In particular, \(\gamma_y^{[m]}(0)=I\) and
\(\gamma_y^{[m]}(mT_y)=M_y^m\).

We recall the convention for the Maslov-type index
from~\cite{LongBott,LongBook}. Let
\(\gamma:[0,T]\to\operatorname{Sp}(2k)\) start at the identity and
write \(M=\gamma(T)\). For \(\omega\in S^1\) set
\[
 \mathcal D_\omega=\{A\in\operatorname{Sp}(2k):\det(A-\omega I)=0\},
 \qquad \nu_\omega(\gamma)=\dim_{\mathbb C}\ker_{\mathbb C}(M-\omega I).
\]
Here \(\operatorname{Sp}(2k)\) is the real symplectic group. If
\(\nu_\omega(\gamma)=0\), the index \(i_\omega(\gamma)\) is the
oriented intersection number with \(\mathcal D_\omega\) of the path
obtained by preceding \(\gamma\) with
\(s\mapsto\operatorname{diag}((2-s)I_k,(2-s)^{-1}I_k)\),
\(0\le s\le1\). At a regular point \(A\) the positive coorientation
is the direction of \(Ae^{sJ_0}\) for increasing \(s\) at zero,
where \(J_0(u,v)=(-v,u)\). If \(\nu_\omega(\gamma)>0\), use the lower
extension
\[
 i_\omega(\gamma)=\sup_{\mathcal U\ni\gamma}
 \inf\{i_\omega(\beta):\beta\in\mathcal U,\ \nu_\omega(\beta)=0\},
\]
where the neighborhoods are taken in the uniform topology on paths
starting at \(I\). See~\cite[Ch.~5]{LongBook} and
\cite[Definition~3.1]{DF}.
In these conventions \(i(y,m)=i_1(\gamma_y^{[m]})\).

The transverse index \(\mu_y(m)\) in \eqref{eq:grading} is the
Conley--Zehnder index of the path on \(\xi\) in a frame extending
over a capping disk for \(y^m\) in \(\Sigma\). Such disks exist and
give the same index since \(\Sigma\) is diffeomorphic to \(S^{2n-1}\)
and \(\pi_1(S^{2n-1})=\pi_2(S^{2n-1})=0\).
The transverse endpoint is nondegenerate by \textup{(H1)}, so this
index equals \(i_1\) of the transverse path.

\subsection{Splitting numbers and the iteration formula}\label{subsec:splitting}

For a path with endpoint \(M\), define the splitting numbers by
\[
 S_M^\pm(\omega)=\lim_{\varepsilon\downarrow0}
   \bigl(i_{\omega e^{\pm i\varepsilon}}(\gamma)-i_\omega(\gamma)\bigr).
\]
They depend only on \(M\) and \(\omega\), vanish off the spectrum,
and add under symplectic direct sums. They are constant on
\(\Omega^0(M)\), the path component of \(M\) among symplectic matrices
with the same eigenvalues on \(S^1\) and the same complex geometric
multiplicities at these eigenvalues~\cite[\S4]{LongBott}.

We compute the contribution of the longitudinal plane
\(E=\operatorname{span}\{\partial_r,R_\lambda\}\).
At \(r=1\), \(\omega_0(\partial_r,R_\lambda)=1\) and
\(E^{\omega_0}=\xi\). Differentiating \(\phi_H^t\) gives
\(d\phi_H^t(\partial_r)=\partial_r+t h''(1)R_\lambda\) and
\(d\phi_H^t(R_\lambda)=R_\lambda\), with the vectors on the right
evaluated at the image point. In the moving symplectic frame
\((R_\lambda,-\partial_r)\), the path on \(E\) is the shear
\[
 \eta(t)=N_1(1,-t h''(1)),\quad \text{ where }
 N_1(1,b)=\begin{pmatrix}1&b\\0&1\end{pmatrix}.
\]
At every \(mT_y\) the shear parameter \(b\) is positive because
\(1<\alpha<2\). The splitting numbers
\(S^+_{N_1(1,b)}(1)=S^-_{N_1(1,b)}(1)=1\) for
\(b>0\)~\cite[List~9.1.12]{LongBook}. For nonreal \(\omega\) near \(1\),
the reference path and \(\eta|_{[0,mT_y]}\) avoid \(\mathcal D_\omega\).
Thus \(i_\omega(\eta|_{[0,mT_y]})=0\) and the definition above gives
\(i_1(\eta|_{[0,mT_y]})=-1\). A capping-disk frame of \(\xi\),
completed by the frame of \(E\), extends over the disk and hence is
homotopic on the boundary to the ambient trivialization. Additivity of
\(i_1\) now yields \(i(y,m)=\mu_y(m)-1\). At the endpoint \(M_y\)
the transverse part \(P_y\) has no eigenvalue \(1\), so
\(S_{M_y}^+(1)=1\).

To estimate the contributions at other unit eigenvalues, deform
\(P_y\) within \(\Omega^0(P_y)\) to its basic normal
form~\cite[Theorem~1.8.10]{LongBook}. Under \textup{(H1)} the possible blocks
are
\[
 D(a)=\begin{pmatrix}a&0\\0&a^{-1}\end{pmatrix},\quad
 R(\theta)=\begin{pmatrix}\cos\theta&-\sin\theta\\
                          \sin\theta&\cos\theta\end{pmatrix},\quad
 N_2(e^{i\theta},B)=\begin{pmatrix}R(\theta)&B\\0&R(\theta)\end{pmatrix}.
\]
Here \(a=\pm2\), \(0<\theta<2\pi\), and \(\theta/(2\pi)\) is
irrational. In the last block \(B=(b_{ij})\) is a real two-by-two
matrix with \(R(\theta)^{\mathsf T}B\) symmetric and \(b_{12}\ne b_{21}\).
The block is nontrivial if \((b_{12}-b_{21})\sin\theta<0\) and trivial
if this quantity is positive.
A rotation has negative splitting number one at \(e^{i\theta}\)
and zero at \(e^{-i\theta}\). A nontrivial \(N_2\) contributes one
at each of these two eigenvalues, and a trivial \(N_2\) contributes
zero at both. The blocks \(D(a)\) have no eigenvalues on \(S^1\).
See~\cite[Lemma~3.2 and Theorem~3.3]{DF}.

Enumerate the distinct \(\theta\in(0,2\pi)\) with
\(S^-_{M_y}(e^{i\theta})>0\) as \(2\pi\alpha_{yj}\),
\(1\le j\le t_y\), and put
\[
 w_{yj}=S^-_{M_y}(e^{2\pi i\alpha_{yj}}),\qquad
 \kappa_y=\sum_{j=1}^{t_y} w_{yj},\qquad
 A_y=i(y,1)+1-\kappa_y .
\]
Each transverse block contributes at most half its dimension to
\(\kappa_y\), so \(\kappa_y\le(2n-2)/2=n-1\).
The splitting-number form of the precise iteration
formula~\cite[Theorem~2.1]{LZ}, also recorded
in~\cite[Theorem~9.3.1]{LongBook}, gives
\[
 i(y,m)=mA_y+2\sum_j w_{yj}\lceil m\alpha_{yj}\rceil-1-\kappa_y .
\]
Here and below the angular sums for \(y\) are over \(1\le j\le t_y\),
and are zero when \(t_y=0\).
Because each \(\alpha_{yj}\) is irrational,
\(\lceil m\alpha_{yj}\rceil=\lfloor m\alpha_{yj}\rfloor+1\).
Adding one to the expression for \(i(y,m)\) gives the first identity
below. Dividing by \(m\) and passing to the limit gives the second.
\begin{align}
 \mu_y(m)&=mA_y+2\sum_j w_{yj}\lfloor m\alpha_{yj}\rfloor+\kappa_y,
       \label{eq:iteration}\\
 \Delta_y&=A_y+2\sum_j w_{yj}\alpha_{yj}. \label{eq:mean}
\end{align}
Consequently,
\[
 \mu_y(m)-m\Delta_y
   =\sum_jw_{yj}\bigl(1-2\{m\alpha_{yj}\}\bigr).
\]
Since \(|1-2\{m\alpha_{yj}\}|\le1\), we obtain the estimate
\begin{equation}\label{eq:error}
 |\mu_y(m)-m\Delta_y|\le\kappa_y\le n-1.
\end{equation}
The variational Morse degree satisfies
\begin{equation}\label{eq:morse-degree}
 I(y^m)=i(y,m)-n=\mu_y(m)-(n+1)
\end{equation}
by~\cite[Theorem~3.10]{DF}. In the counting arguments below the
\emph{degree} of \(y^m\) is \(\mu_y(m)\). Its sign in the
variational Euler characteristic is
\((-1)^{I(y^m)}=(-1)^{\mu_y(m)-(n+1)}\).

\subsection{Morse coefficients and orbit types}\label{subsec:morse}

Recall that \(G_y\) is the set of good iteration numbers of \(y\).
By \eqref{eq:iteration},
\begin{align}\label{sec2:1}
	 \mu_y(m)-\mu_y(1)\equiv(m-1)A_y\pmod 2,\qquad
	G_y=
	\begin{cases}
		\N,& A_y\ \text{even},\\
		\{1,3,5,\ldots\},& A_y\ \text{odd}.
	\end{cases}
\end{align}
In particular, for even \(k\) and \(1\le\ell<k\), membership in \(G_y\)
is the same for \(\ell\), \(k-\ell\), and \(k+\ell\).
Define the Morse coefficient in shifted degree \(j\in\Z\) by
\begin{align}\label{sec2:3}
	 C_j=\#\{(y,m):y\in\PP,\ m\in G_y,\ \mu_y(m)=j\}.
\end{align}
The pair \((y,m)\) specifies a prime orbit and a covering multiplicity.
It is counted once, independently of the choice of starting point on
the orbit. By \eqref{eq:error}, \(\mu_y(m)=j\) implies
\(m\delta_*\le j+n-1\). Thus \(C_j\) is finite and vanishes for
\(j\le-(n-1)\).

The circle acts on parametrized loops by time translation. The
translates of \(y^m\) form a single critical orbit with cyclic
stabilizer of order \(m\). With rational coefficients its local
equivariant critical module, in the normalized variational grading,
has rank one in degree \(I(y^m)\) when \(m\in G_y\). All its other
ranks are zero, and a bad iterate contributes zero in every degree.
This follows from the local module calculation
in~\cite[Proposition~4.2 and Remark~4.4]{LLW}, with the index
normalization in \eqref{eq:morse-degree}. In the
variational grading the Betti series is
\(1/(1-t^2)=1+t^2+t^4+\cdots\). The shift
\(\mu_y(m)=I(y^m)+n+1\) therefore gives Betti series
\(t^{n+1}/(1-t^2)\). In this grading the Morse relation reads
\begin{equation}\label{eq:morse}
 C_j=B_j+U_j+U_{j-1},\qquad U_j\in\Z_{\ge0},
 \qquad
 B_j=\begin{cases}
  1,&j\in n+1+2\Z_{\ge0},\\
  0,&\text{otherwise}.
 \end{cases}
\end{equation}
Here \(B_j\) is the Betti number in shifted degree \(j\). The
nonnegative integers \(U_j\) are the coefficients of the remainder
series in the equivariant Morse identity
\[
 \sum_{j\in\Z}C_jt^j
 =\frac{t^{n+1}}{1-t^2}+(1+t)\sum_{j\in\Z}U_jt^j.
\]
Here \(t\) is a formal variable. The underlying Morse relation
is proved in~\cite[Corollary~5.11 and (8.12)]{Viterbo}.
The finite-cutoff relation and its limit, when all mean indices
are nonzero, are given in~\cite[(6.6), (6.18)]{LLW}.
The formula above shifts those degrees by \(n+1\).

For completeness, the coefficients in \eqref{eq:morse} can be
chosen as one globally defined sequence. At a finite period cutoff,
the Morse relation has an additional term \(Q_a(t)\), as in
\cite[(6.6)]{LLW}. On each fixed finite degree interval, its
coefficients vanish when the cutoff is sufficiently large
\cite[Theorem~5.1 and (6.4)]{LLW}. Indeed,
\eqref{eq:error} and \eqref{eq:morse-degree} give
\(I(y^m)\ge m\delta_*-2n\), and there are only finitely many
prime orbits. Iterates beyond a sufficiently large period cutoff
therefore have indices above the prescribed interval.
The coefficients counting critical modules also stabilize on
that interval. This gives the limiting relation with nonnegative
integer coefficients.

Since the terms on the right of \eqref{eq:morse} are nonnegative,
\(C_j=0\) implies \(U_j=U_{j-1}=0\). In particular \(U_j=0\) for
\(j\le-(n-1)\). Put \(j_*=-(n-1)\). Recursion in
\eqref{eq:morse} gives
\begin{equation}\label{eq:fixed-remainder}
 U_j=\sum_{\ell=j_*+1}^{j}(-1)^{j-\ell}(C_\ell-B_\ell)
 \qquad(j>j_*).
\end{equation}
Thus \(U_j\) is uniquely determined by the global coefficients
\(C_\ell\). It does not depend on a period cutoff or on any
common index jump chosen later. This independence allows us to
choose a minimum of the relevant coefficients before selecting the
jump data for that degree in Lemma~\ref{lem:window}.

For use in Section~\ref{sec:smaller}, note that \textup{(H3)}
is equivalent to
\begin{equation}\label{eq:zero}
 C_0=0.
\end{equation}

Write \(\sigma_y=(-1)^{\mu_y(1)-(n+1)}\), the sign of each good
iterate in the variational grading. The mean Euler characteristic
of \(y\) is
\[
 \mchi(y)=\sigma_y\lim_{r\to\infty}
              \frac{\#(G_y\cap[1,r])}{r}
 \ \in\ \{1,\tfrac12,-1,-\tfrac12\}.
\]
Indeed, \(\#(G_y\cap[1,r])/r\) tends to \(1\) or \(1/2\)
according as \(A_y\) is even or odd. Multiplication by \(\sigma_y\)
gives \(\mchi(y)\), in agreement with~\cite[Remark~4.9]{LLW}.
Let
\[
 \PP_\eta=\{y\in\PP:\mchi(y)=\eta\},\qquad
 (a,b,c,h)=
 \bigl(\#\PP_1,\#\PP_{1/2},\#\PP_{-1},\#\PP_{-1/2}\bigr).
\]
The four sets partition \(\PP\), so \(q=a+b+c+h\). Their contributions
to the coefficient counts are described by
\[
\begin{array}{c|c|c|c}
\text{class} & \text{number of prime orbits} & G_y
             & \mu_y(m)\ (m\in G_y)\\ \hline
\PP_1       & a & \N               & \equiv n+1\pmod2\\
\PP_{1/2}   & b & \{1,3,5,\ldots\} & \equiv n+1\pmod2\\
\PP_{-1}    & c & \N               & \equiv n\pmod2\\
\PP_{-1/2}  & h & \{1,3,5,\ldots\} & \equiv n\pmod2
\end{array}
\]
The parity condition on the iteration number is independent of
the parity of \(n\). Only the degree parity in the last column
depends on \(n\).
When \(k\) is even, exactly \(k\) or \(k/2\) of the first \(k\)
iterates are good. Hence
\begin{equation}\label{eq:complete-count}
 \sum_{\substack{1\le m\le k\\y^m\ {\rm good}}}
 (-1)^{\mu_y(m)-(n+1)}=k\mchi(y).
\end{equation}
\Needspace{6\baselineskip}
We record also the resonance identity of Liu, Long, and
Wang~\cite[Theorem~1.1]{LLW}, which will be used only in
Section~\ref{sec:smaller}. Under the present positivity assumption,
it reads
\begin{equation}\label{eq:resonance}
 \sum_{y\in\PP}\frac{\mchi(y)}{\Delta_y}=\frac12.
\end{equation}
With the variational sign in \eqref{sec2:2}, the right-hand side
is \(1/2\) for both parities of \(n\).

\section{A common index jump}\label{sec:jump}

In this section we specialize the generalized common index jump theorem~\cite[Theorem~3.6]{DLLW} to nondegenerate closed characteristics with positive mean indices in the grading 
\(\mu_y(m)=i(y,m)+1\).
We include a direct proof to specify the approximation errors
and integer identities used in the construction of complementary
normalized jump data in Section~\ref{sec:smaller}.

For \(L\in\N\), \emph{jump data of order \(L\)} are a family
\[
 J=\bigl(d,(k_y,s_y)_{y\in\PP}\bigr)
\]
with \(d\) a positive even integer, \(k_y>L\) positive and even,
and \(s_y=\mu_y(k_y)-d\), satisfying
\eqref{eq:jump}--\eqref{eq:offset} below.
The integer \(d\) is a common degree, called the \emph{center}
of the data, whereas \(k_y\) depends on \(y\).
Both \(d\) and \(k_y\) are even for every \(n\ge2\).
Relative to these data, \(y^{k_y}\) is the
\emph{central iterate}. Every \(y^m\) with \(m\ne k_y\) is
\emph{noncentral}. For \(1\le\ell\le L\), we compare the
reference iterate \(y^\ell\) with \(y^{k_y-\ell}\) and
\(y^{k_y+\ell}\). These latter iteration numbers lie to the
left and right of \(k_y\); the order of their indices can differ
from the order of the iteration numbers when \(\mu_y(\ell)<0\).

\begin{lemma}[{cf. \cite[Theorem~3.6]{DLLW}}]\label{lem:recurrence}
Assume \textup{(H1)}, finiteness of \(\PP\), and positivity of
all mean indices. For every \(L,H\in\N\) there are jump data
\(J=(d,(k_y,s_y)_{y\in\PP})\) of order \(L\), with
\(d=2N\ge H\) and \(k_y>L\) even, such that
\begin{align}
 \mu_y(k_y\pm\ell)&=d\pm\mu_y(\ell)
       &&(y\in\PP,\ 1\le\ell\le L),\label{eq:jump}\\
 \mu_y(k_y)&=d+s_y,\qquad |s_y|\le n-1
       &&(y\in\PP).\label{eq:offset}
\end{align}
\end{lemma}

\begin{proof}
\emph{Choice of the integers.}
Fix an ordering of the finite set \(\PP\) and form the vector
\[
 v=\left((1/\Delta_y)_{y\in\PP},
         (2\alpha_{yj}/\Delta_y)_{y\in\PP,\,1\le j\le t_y}\right)
 \in\R^M,\qquad M=q+\sum_{y\in\PP}t_y.
\]
For every \(\varepsilon>0\) and every integer \(T\ge1\), there is an
integer \(N\ge T\) such that
\[
 \max_{1\le i\le M}\|Nv_i\|_{\Z}<\varepsilon.
\]
Indeed, take an integer \(Q>1/\varepsilon\) and partition
\([0,1)^M\) into \(Q^M\) half-open cubes of side \(1/Q\). Two of the
\(Q^M+1\) points \((\{rTv_i\})_{i=1}^M\), \(0\le r\le Q^M\), lie
in the same cube. Let \(r_1<r_2\) be the corresponding values
of \(r\). Then \(N=T(r_2-r_1)\ge T\), and
\[
 \|Nv_i\|_{\Z}\le|\{r_2Tv_i\}-\{r_1Tv_i\}|
       <1/Q<\varepsilon\qquad(1\le i\le M).
\]

Choose \(0<\varepsilon<1/6\) such that
\begin{align*}
 \varepsilon(\Delta_y+3\kappa_y)&<1 &&(y\in\PP),\\
 6\varepsilon&<\|\ell\alpha_{yj}\|_{\Z}
       &&(y\in\PP,\ 1\le j\le t_y,\ 1\le\ell\le L).
\end{align*}
Every distance in the last line is positive by irrationality, and
there are only finitely many conditions. Pick \(T\in\N\) with
\(2T\ge H\) and \(2T/\Delta_y-1>L\) for every \(y\).
For the integer \(N\ge T\) supplied by the approximation above, set
\[
 d=2N,\qquad k_y=2\operatorname{rnd}(N/\Delta_y),\qquad
 b_{yj}=\operatorname{rnd}(2N\alpha_{yj}/\Delta_y).
\]
Then \(d\ge H\) and \(k_y\ge2N/\Delta_y-1>L\).
The errors \(u_y=N/\Delta_y-k_y/2\) and
\(\beta_{yj}=2N\alpha_{yj}/\Delta_y-b_{yj}\) satisfy
\(|u_y|,|\beta_{yj}|<\varepsilon\). It follows that
\[
 e_{yj}:=k_y\alpha_{yj}-b_{yj}
       =\beta_{yj}-2\alpha_{yj}u_y,\qquad |e_{yj}|<3\varepsilon.
\]

\smallskip
\noindent\emph{The index identities.}
By \eqref{eq:mean},
\[
 k_yA_y+2\sum_jw_{yj}b_{yj}-d
      =-2\Delta_yu_y-2\sum_jw_{yj}e_{yj}.
\]
The left side is an even integer, and the right side has absolute
value less than \(2\varepsilon(\Delta_y+3\kappa_y)<2\). Hence
\begin{equation}\label{eq:integer-center}
 k_yA_y+2\sum_jw_{yj}b_{yj}=d.
\end{equation}
For \(1\le\ell\le L\), addition of \(e_{yj}\) to
\(\pm\ell\alpha_{yj}\) leaves its integer part unchanged, because
\(|e_{yj}|<\|\ell\alpha_{yj}\|_{\Z}\). Thus
\begin{align*}
 \lfloor(k_y+\ell)\alpha_{yj}\rfloor
    &=b_{yj}+\lfloor\ell\alpha_{yj}\rfloor,\\
 \lfloor(k_y-\ell)\alpha_{yj}\rfloor
    &=b_{yj}-\lfloor\ell\alpha_{yj}\rfloor-1 .
\end{align*}
Substituting the first formula into \eqref{eq:iteration} gives
\(d+\ell A_y+2\sum_jw_{yj}\lfloor\ell\alpha_{yj}\rfloor+\kappa_y
=d+\mu_y(\ell)\). In the second formula the terms \(-1\) give
\(-2\sum_jw_{yj}=-2\kappa_y\), so that
\(\mu_y(k_y-\ell)=d-\mu_y(\ell)\).

For the central iterate, \(e_{yj}\ne0\) since \(k_y>0\) and
\(\alpha_{yj}\) is irrational. As \(|e_{yj}|<1\), we have
\(\lfloor k_y\alpha_{yj}\rfloor=b_{yj}-\mathbf{1}_{\{e_{yj}<0\}}\).
It follows from \eqref{eq:iteration} and \eqref{eq:integer-center} that
\[
 s_y=\mu_y(k_y)-d
     =\kappa_y-2\sum_{j:\,e_{yj}<0}w_{yj}
     =\sum_jw_{yj}\operatorname{sgn}(e_{yj}),\qquad
 |s_y|\le\kappa_y\le n-1 .
\]
This proves the assertion.
\end{proof}

The index estimate \eqref{eq:error} at \(m=k_y\) gives
\[
 |k_y\Delta_y-d|\le|k_y\Delta_y-\mu_y(k_y)|+|s_y|\le2(n-1).
\]
For \(t\in\Z\) with \(k_y+t\ge1\), applying \eqref{eq:error} at
\(m=k_y+t\) consequently gives
\begin{equation}\label{eq:far}
 |\mu_y(k_y+t)-d-t\Delta_y|\le 3(n-1).
\end{equation}
In particular, if \(|t|>3(n-1)/\delta_*\), then
\(\mu_y(k_y+t)-d\) has the sign of \(t\), since
\(|t|\Delta_y\ge|t|\delta_*>3(n-1)\).

\section{Morse coefficient bounds and asymptotic counting}
\label{sec:counting}

We first retain the boundary terms when comparing alternating
Morse sums over two finite intervals related by a common index jump.
The resulting identity, applied at a minimum of the global Morse
remainder, gives a uniform lower bound for \(C_p\) in every degree
\(p\ge n\) with \(p\equiv n\pmod2\). Summing these bounds gives
a lower bound on the asymptotic density of good iterates in those
degrees. We then compute the same density by counting iterates
separately on their prime orbits and obtain a quantitative bound
for the number of prime characteristics.

Recall that \(C_j\) counts good iterates of degree \(\mu_y(m)=j\),
and \(U_j\) is the globally defined nonnegative remainder in
\eqref{eq:morse}. The numbers \(a,b,c,h\) count the classes
\(\PP_1,\PP_{1/2},\PP_{-1},\PP_{-1/2}\), respectively. A jump datum
\(J\) consists of \(d\) and the pairs \((k_y,s_y)\) for all
\(y\in\PP\), as in Section~\ref{sec:jump}.

\subsection{A minimum argument for the global Morse remainder}

\begin{lemma}\label{lem:window}
Let \(n\ge2\). Suppose that \textup{(H1)} holds, that
\(\PP\) is finite, and that \(\Delta_y>0\) for every \(y\in\PP\). Set
\begin{equation}\label{eq:DK}
 D=n-a+c.
\end{equation}
For every integer \(p\ge n\) with \(p\equiv n\pmod2\),
\begin{equation}\label{eq:coefficient-lower}
 C_p\ge U_{p-1}\ge\max\{0,D\}.
\end{equation}
\end{lemma}

\begin{proof}
Fix an integer \(p\ge n\) with \(p\equiv n\pmod2\), and choose an integer
\(L>(p+3(n-1))/\delta_*\).
Choose jump data \(J=(d,(k_y,s_y)_{y\in\PP})\) of order \(L\) from
Lemma~\ref{lem:recurrence}, with \(d-p+1\ge n+1\). For \(j\in\Z\) put
\[
 Z_j=Z_j(J)=\#\{y\in\PP:y^{k_y}\text{ is good and }s_y=j\}.
\]
We first count the good iterates in degree \(d+j\). For \(|j|<p\),
\begin{equation}\label{eq:reflection}
 C_{d+j}=C_j+C_{-j}+Z_j .
\end{equation}
Indeed, if \(\mu_y(m)=d+j\) and \(t=m-k_y\), then \eqref{eq:far} gives
\[
 |t|\delta_*\le |t|\Delta_y
       \le |j|+3(n-1)<p+3(n-1).
\]
Hence \(|t|<L\). The case \(t=0\) accounts for \(Z_j\) when we count
good iterates. If \(t>0\), \eqref{eq:jump} gives
\(\mu_y(t)=j\). If \(t<0\), the same formula gives \(\mu_y(-t)=-j\).
Conversely, if \(\mu_y(\ell)\in\{j,-j\}\) for some \(\ell\in\N\), then
\(\ell\delta_*\le |j|+n-1<L\delta_*\) by \eqref{eq:error}.
In particular \(k_y-\ell\ge1\). The two maps
\[
 (y,\ell)\longmapsto(y,k_y+\ell)
      \quad\text{for }\mu_y(\ell)=j,\qquad
 (y,\ell)\longmapsto(y,k_y-\ell)
      \quad\text{for }\mu_y(\ell)=-j
\]
therefore give bijections onto the pairs in degree \(d+j\) with
\(m>k_y\) and \(m<k_y\), respectively. Both maps preserve goodness,
since \(\mu_y(k_y\pm\ell)=d\pm\mu_y(\ell)\) and \(d\) is even.
The two images are disjoint also when \(j=0\), which proves
\eqref{eq:reflection} in that degree as well.

Put \(\varepsilon_n=(-1)^{n+1}\). For an even integer \(r\), define
\[
 W_p(r)=\sum_{j=1-p}^{p-1}(-1)^j C_{r+j}.
\]
We compare the two window sums \(W_p(0)\) and \(W_p(d)\).
Every offset satisfies \(|s_y|\le n-1<p\). At the even iteration
number \(k_y\), the good iterates come precisely from
\(\PP_1\cup\PP_{-1}\). On \(\PP_1\), their degree \(d+s_y\)
is congruent to \(n+1\) modulo two; on \(\PP_{-1}\), it is
congruent to \(n\). Since \(d\) is even, the central contribution is
\[
 \sum_{|j|<p}(-1)^j Z_j=\varepsilon_n(a-c).
\]
The two sums of \(C_j\) and \(C_{-j}\) in \eqref{eq:reflection}
agree after changing \(j\) to \(-j\), so
\begin{equation}\label{eq:signed-reflection}
 W_p(d)=2W_p(0)+\varepsilon_n(a-c).
\end{equation}

We next express these sums using the Morse relation, retaining
both boundary terms. For any integers \(u\le v\), reindexing the sum of
\((-1)^jU_{j-1}\) in \eqref{eq:morse} gives
\[
 \sum_{j=u}^{v}(-1)^j C_j
   =\sum_{j=u}^{v}(-1)^j B_j
     +(-1)^vU_v+(-1)^uU_{u-1}.
\]
In the interval \([1-p,p-1]\), the number of Betti degrees is
\((p-n)/2\). It is zero for \(p=n\), and for \(p>n\) the Betti
degrees are \(n+1,n+3,\ldots,p-1\). The interval
\([d-p+1,d+p-1]\) contains \(p\) degrees congruent to \(n+1\)
modulo two, all at least \(n+1\). Each nonzero Betti term
therefore has sign \(\varepsilon_n\). The endpoints of both
intervals also have this sign, because \(p\equiv n\pmod2\)
and \(d\) is even. Finally, \(U_{-p}=0\), since \(p\ge n\).
It follows that
\begin{align*}
 W_p(0)&=\varepsilon_n\left(\frac{p-n}{2}+U_{p-1}\right),\\
 W_p(d)&=\varepsilon_n\left(p+U_{d+p-1}+U_{d-p}\right).
\end{align*}
Substituting these expressions in \eqref{eq:signed-reflection}
and dividing by \(\varepsilon_n\) gives
\begin{equation}\label{eq:boundary}
 U_{d+p-1}+U_{d-p}=2U_{p-1}+a-c-n.
\end{equation}
The identity holds for every \(p\ge n\) with
\(p\equiv n\pmod2\), after choosing the jump data for that
value of \(p\).

Nonnegativity of the left-hand side of \eqref{eq:boundary} gives
\(U_{p-1}\ge D/2\). We now use a minimum of the global remainder
to obtain the bound \(U_{p-1}\ge D\). The coefficients \(U_j\)
form the fixed sequence described in \eqref{eq:fixed-remainder}.
The nonempty set
\(\{U_{p-1}:p\ge n,\ p\equiv n\pmod2\}\) consists of nonnegative
integers, so it has a minimum. Set
\[
 e=\min_{\substack{p\ge n\\p\equiv n\ (\mathrm{mod}\ 2)}}U_{p-1},
 \qquad U_{p_0-1}=e
\]
for some \(p_0\ge n\) with \(p_0\equiv n\pmod2\). Fix this
\(p_0\), and apply \eqref{eq:boundary} with the corresponding
even center \(d\). The integer \(d+p_0\) is at least \(n\)
and is congruent to \(n\) modulo two, so \(U_{d+p_0-1}\ge e\).
For the remaining term \(U_{d-p_0}\), whose index has the opposite
parity, we use only nonnegativity. We therefore obtain
\[
 e\le U_{d+p_0-1}
   \le U_{d+p_0-1}+U_{d-p_0}=2e-D.
\]
It follows that \(e\ge D\). Thus \(U_{p-1}\ge\max\{0,D\}\)
for every \(p\ge n\) with \(p\equiv n\pmod2\). For these
degrees \(B_p=0\), and
\eqref{eq:morse} gives \(C_p=U_p+U_{p-1}\ge U_{p-1}\).
\end{proof}

\begin{unnumremark}
The coefficient \(C_0\) need not vanish. At \(j=0\), the
reflection identity is \(C_d=2C_0+Z_0(J)\): a good reference
iterate of degree zero gives two distinct iterates at degree
\(d\). Both contributions are retained in the alternating sum,
so the proof of Lemma~\ref{lem:window} does not require \(C_0=0\)
or any other individual index exclusion. Once a minimizing degree
\(p_0\) has been chosen, only one set of jump data is needed.
\end{unnumremark}

\begin{unnumremark}
When \(n=2\), the smallest allowed value is \(p=2\). The low
interval \([-1,1]\) contains no Betti degree, whereas
\([d-1,d+1]\), for even \(d\ge4\), contains the two Betti
degrees \(d-1\) and \(d+1\). Both have negative sign in the
alternating sum. Thus the same calculation applies, and the
minimum is taken over the coefficients \(U_j\) with odd
\(j\ge1\).
\end{unnumremark}

\subsection{Counting good iterates in the selected degrees}

For real \(T\ge n\), define
\begin{equation}\label{eq:rho}
 \begin{split}
 A(T)&=\sum_{\substack{n\le p\le T\\p\equiv n\ (\mathrm{mod}\ 2)}}C_p,\\
 \rho&=\sum_{y\in\PP_{-1}}\frac1{\Delta_y}
       +\sum_{y\in\PP_{-1/2}}\frac1{2\Delta_y}.
 \end{split}
\end{equation}
Thus \(A(T)\) counts pairs \((y,m)\) representing good iterates
whose indices lie in \([n,T]\cap(n+2\Z)\). These are the
degrees in which the Betti coefficients vanish. Different iterates
of the same degree are counted separately. The next lemma
identifies \(\rho\) as the coefficient of \(T\) in the
asymptotic formula for \(A(T)\) and turns the coefficient bound
into a lower bound for this density.

\Needspace{11\baselineskip}
\begin{lemma}\label{lem:counting}
Under the hypotheses of Lemma~\ref{lem:window}, one has
\begin{equation}\label{eq:asymptotic}
 A(T)=\rho T+O(1)\qquad(T\to\infty).
\end{equation}
In particular,
\begin{equation}\label{eq:rho-lower}
 2\rho\ge D.
\end{equation}
If \(\delta>0\) satisfies \(\Delta_y\ge\delta\) for every
\(y\in\PP_{-1}\cup\PP_{-1/2}\), then
\begin{equation}\label{eq:rho-upper}
 2\rho\le\frac{2c+h}{\delta}.
\end{equation}
\end{lemma}

\begin{proof}
The prime orbits whose good iterates have degrees congruent to
\(n\) modulo two form the set
\(\mathcal O=\PP_{-1}\cup\PP_{-1/2}\).
For \(y\in\mathcal O\), put \(g_y=|\mchi(y)|\). Thus \(g_y=1\)
on \(\PP_{-1}\) and \(g_y=1/2\) on \(\PP_{-1/2}\).
For real \(x\ge0\), the number of good iteration numbers up to
\(x\) is
\[
 F_y(x)=\#\{m\in G_y:m\le x\}
 =\begin{cases}
   \lfloor x\rfloor,&y\in\PP_{-1},\\
   \lfloor(x+1)/2\rfloor,&y\in\PP_{-1/2}.
 \end{cases}
\]
Indeed, every positive iteration number is good in the first case,
and precisely the odd iteration numbers are good in the second,
by \eqref{sec2:1} and \eqref{sec2:2}. In both cases,
\(\lvert F_y(x)-g_yx\rvert\le1\).

For \(u\in\R\), write
\(N_y(u)=\#\{m\in G_y:\mu_y(m)\le u\}\).
This number is finite, since \eqref{eq:error} gives
\(m\Delta_y\le u+n-1\) for each exponent counted by \(N_y(u)\).
For \(T\ge n\), an exponent with
\(m\Delta_y\le T-(n-1)\) satisfies \(\mu_y(m)\le T\), whereas
every exponent with \(\mu_y(m)\le T\) satisfies
\(m\Delta_y\le T+(n-1)\), again by \eqref{eq:error}.
Consequently,
\[
 F_y\!\left(\frac{T-(n-1)}{\Delta_y}\right)
 \le N_y(T)\le
 F_y\!\left(\frac{T+(n-1)}{\Delta_y}\right).
\]
Combining these inequalities with the bound for \(F_y\) yields
\[
 \left|N_y(T)-\frac{g_y}{\Delta_y}T\right|
 \le 1+\frac{(n-1)g_y}{\Delta_y}.
\]
The good iterates of \(y\) below degree \(n\) account for exactly
\(N_y(n-1)\) terms. Each satisfies \(m\Delta_y\le2(n-1)\), so
\(N_y(n-1)\le\lfloor2(n-1)/\Delta_y\rfloor\).
All good degrees of \(y\) are congruent to \(n\) modulo two.
Consequently its contribution to \(A(T)\) is
\[
 Q_y(T)=\#\{m\in G_y:n\le\mu_y(m)\le T\}
       =N_y(T)-N_y(n-1).
\]
It follows that
\[
 \left|Q_y(T)-\frac{g_y}{\Delta_y}T\right|
 \le 1+\frac{(n-1)g_y}{\Delta_y}+N_y(n-1)
 \qquad(T\ge n).
\]
Now \(A(T)=\sum_{y\in\mathcal O}Q_y(T)\), and the two values of
\(g_y\) give \(\sum_{y\in\mathcal O}g_y/\Delta_y=\rho\).
The set \(\mathcal O\) is finite, so summing the last inequality
gives a bound for \(\lvert A(T)-\rho T\rvert\) independent of
\(T\). This proves \eqref{eq:asymptotic}, and division by \(T\)
gives \(\lim_{T\to\infty}A(T)/T=\rho\).

Let \(\nu_n(T)=1+\lfloor(T-n)/2\rfloor\) be the number of
integers in \([n,T]\cap(n+2\Z)\). Lemma~\ref{lem:window} gives \(C_p\ge D\)
for each such integer, hence \(A(T)\ge D\nu_n(T)\).
Dividing by \(T\) and using \(\nu_n(T)/T\to1/2\) gives
\(\rho\ge D/2\), which proves \eqref{eq:rho-lower}.
If \(\Delta_y\ge\delta>0\) on \(\mathcal O\), substituting
\(1/\Delta_y\le1/\delta\) in \eqref{eq:rho} gives
\(\rho\le c/\delta+h/(2\delta)\), proving
\eqref{eq:rho-upper}.

\end{proof}

\begin{corollary}\label{cor:quantitative}
Under the hypotheses of Lemma~\ref{lem:window}, write
\(q=\#\PP\) and \(b=\#\PP_{1/2}\). Then
\begin{equation}\label{eq:quantitative-multiplicity}
 q\ge n+b+
 2\sum_{\substack{y\in\PP\\\mchi(y)<0}}
 |\mchi(y)|\left(1-\frac1{\Delta_y}\right).
\end{equation}
In particular, if \textup{(H2)} holds, then \(q\ge n+b\).
\end{corollary}

\begin{proof}
Since \(q=a+b+c+h\) and \(D=n-a+c\),
\eqref{eq:rho-lower} gives
\[
 q=n+b+2c+h-D\ge n+b+2c+h-2\rho.
\]
The definitions of the orbit classes and \(\rho\) give
\[
 2c+h=2\sum_{\substack{y\in\PP\\\mchi(y)<0}}|\mchi(y)|,
 \qquad
 \rho=\sum_{\substack{y\in\PP\\\mchi(y)<0}}
       \frac{|\mchi(y)|}{\Delta_y}.
\]
Substitution proves \eqref{eq:quantitative-multiplicity}.
Under \textup{(H2)}, each summand in that inequality is
nonnegative, which gives \(q\ge n+b\).
\end{proof}
\begin{corollary}[Rigidity at the minimal orbit count]
\label{cor:minimal-rigidity}
Under the hypotheses of Lemma~\ref{lem:window}, assume that
\(\#\PP=n\) and that
\[
 \Delta_y>1
 \qquad\text{whenever }\mchi(y)<0.
\]
Then \(\mchi(y)=1\) for every \(y\in\PP\), and every iterate
of every prime closed characteristic is good. Moreover,
\[
 U_j=0,\qquad C_j=B_j
 \qquad\text{for every }j\in\Z.
\]
In particular, the Reeb flow is dynamically convex:
\[
 \mu_y(m)\ge n+1
 \qquad\text{for every }y\in\PP\text{ and }m\in\N.
\]
\end{corollary}

\begin{proof}
Applying Corollary~\ref{cor:quantitative} with \(q=n\), we obtain
\[
 0\ge b+
 2\sum_{\substack{y\in\PP\\\mchi(y)<0}}
 |\mchi(y)|\left(1-\frac1{\Delta_y}\right).
\]
Every term on the right is nonnegative. Furthermore, each
characteristic with negative mean Euler characteristic would
contribute a strictly positive term, by the assumed strict
mean index bound. Hence
\[
 b=c=h=0,\qquad a=n.
\]
Thus every prime characteristic belongs to \(\PP_1\).
The orbit classification in Section~\ref{subsec:morse} shows
that every iterate is good and that
\[
 \mu_y(m)\equiv n+1\pmod2
 \qquad(y\in\PP,\ m\in\N).
\]

Consequently, \(C_j=0\) whenever \(j\equiv n\pmod2\).
For these degrees, \(B_j=0\) as well, and the Morse relation
\eqref{eq:morse} gives
\[
 0=C_j=B_j+U_j+U_{j-1}=U_j+U_{j-1}.
\]
Since the remainder coefficients are nonnegative, both
\(U_j\) and \(U_{j-1}\) vanish. As \(j\) ranges over the integers
congruent to \(n\) modulo two, these pairs cover all integer
indices. Therefore \(U_j=0\) for every \(j\in\Z\), and
\eqref{eq:morse} yields \(C_j=B_j\) in every degree.

Finally, every iterate \(y^m\) is good and therefore contributes
to \(C_{\mu_y(m)}\). It follows that
\[
 B_{\mu_y(m)}=C_{\mu_y(m)}\ge1.
\]
By the definition of the Betti coefficients in \eqref{eq:morse},
\[
 \mu_y(m)\in n+1+2\Z_{\ge0}.
\]
In particular, \(\mu_y(m)\ge n+1\) for every iterate, which
proves dynamical convexity.
\end{proof}

\begin{unnumremark}
The strict mean index bound is used to exclude all characteristics
with negative mean Euler characteristic. If it is replaced by
\(\Delta_y\ge1\) on that class, Corollary~\ref{cor:quantitative}
with \(q=n\) gives only
\[
 b=0,\qquad
 \Delta_y=1\quad\text{whenever }\mchi(y)<0.
\]
These conclusions do not exclude the negative mean Euler
characteristic classes. 
\end{unnumremark}
\section{The mean index criterion}\label{sec:main-proof}

\begin{proof}[Proof of Theorem~\ref{thm:main}]
Suppose that \(\Sigma\) has fewer than \(n\) geometrically
distinct prime closed characteristics. Then \(\PP\) is finite,
\(q\le n-1\), and the assumptions of Sections~\ref{sec:prelim}--
\ref{sec:counting} hold. Since \(a\le q\le n-1\),
\[
 D=n-a+c\ge1.
\]
Using \(q=a+b+c+h\), we also have
\begin{equation}\label{eq:negative-count}
 2c+h=q-b-a+c\le n-1-a+c=D-1.
\end{equation}
Lemma~\ref{lem:counting} and \textup{(H2)}, applied only to
\(\PP_{-1}\cup\PP_{-1/2}\), now give
\[
 D\le2\rho\le2c+h\le D-1.
\]
This is a contradiction. Hence \(q\ge n\).
\end{proof}

\begin{proof}[Proof of Corollary~\ref{cor:uniform}]
The uniform bound \(\Delta_y\ge1\) implies \textup{(H2)}.
The conclusion follows from Theorem~\ref{thm:main}.
\end{proof}

\begin{proof}[Proof of Corollary~\ref{cor:absolute}]
If some prime closed characteristic has negative mean index,
\cite[Theorem~B]{CGGFiltered} gives infinitely many prime closed
characteristics. Otherwise \(|\Delta_y|\ge1\) implies
\(\Delta_y\ge1\) for every prime characteristic, and
Corollary~\ref{cor:uniform} applies.
\end{proof}

\begin{remark}
Theorem~\ref{thm:main} and Corollary~\ref{cor:uniform} use the
\(C^3\) variational framework of~\cite{LLW}. In
Corollary~\ref{cor:absolute} we assume smoothness in order to
apply~\cite[Theorem~B]{CGGFiltered}. 
\end{remark}

\begin{remark}[Even half-dimensions with an index exclusion]
\label{rem:even-exclusion}
Suppose that \(n\ge2\) is even, \textup{(H1)} and \textup{(H3)}
hold, and every prime characteristic has positive mean index.
The coefficient estimate then proves the multiplicity bound
without a quantitative mean index assumption. If there are
infinitely many prime characteristics, the conclusion is immediate.
Otherwise choose an integer \(L>(n+3(n-1))/\delta_*\) and jump
data of order \(L\) with even center \(d\ge2n\). The reflection
formula \eqref{eq:reflection}, applied with \(p=n\) and \(j=0\),
gives
\[
 C_d=2C_0+Z_0(J)=Z_0(J),
\]
because \textup{(H3)} gives \(C_0=0\). An orbit counted by
\(Z_0(J)\) has a good even iterate of even degree. Its even
iteration number implies that all its iterates are good, and its
degree parity gives \(\mchi(y)=-1\). Hence \(Z_0(J)\le c\).
Since \(d\ge n\) and \(d\equiv n\pmod2\),
Lemma~\ref{lem:window} yields
\[
 n-a+c=D\le C_d=Z_0(J)\le c.
\]
Thus \(a\ge n\), and in particular \(q\ge n\). This recovers
the positive-mean-index case of the stronger even-dimensional
result in~\cite[Theorem~1.2]{DLLW}, which assumes only that all
mean indices are nonzero.
\end{remark}

\begin{example}[Sharpness of the multiplicity bound]\label{ex:ellipsoid}
Let
\[
 \Sigma=\left\{z\in\mathbb C^n:
             \sum_{j=1}^n\frac{|z_j|^2}{a_j}=1\right\},
 \qquad a_j>0,
\]
where \(a_j/a_k\) is irrational for \(j\ne k\).
The Reeb flow rotates the \(j\)-th coordinate with angular
velocity \(2/a_j\). If a point \(z\) on a periodic orbit of period
\(T\) had two nonzero coordinates \(z_j,z_k\), then both
\(T/(\pi a_j)\) and \(T/(\pi a_k)\) would be positive integers.
Their ratio would be \(a_k/a_j\), contrary to irrationality.
Thus the \(n\) coordinate circles are the only prime closed
characteristics. The \(j\)-th circle \(y_j\) has least period
\(\pi a_j\). Along \(y_j^m\), the transverse return map on the
\(k\)-th coordinate plane, \(k\ne j\), is a rotation through
\(2\pi m a_j/a_k\), whose eigenvalues differ from one.
All iterates are therefore nondegenerate. The rotation indices,
together with the contribution \(2m-1\) from the orbit plane in
the ambient convention, give
\[
 i(y_j,m)=2m-1+
    \sum_{k\ne j}\bigl(2\lfloor ma_j/a_k\rfloor+1\bigr)\ge n.
\]
Moreover,
\[
 \Delta_{y_j}=2\sum_{k=1}^n\frac{a_j}{a_k}>2.
\]
Thus these hypersurfaces satisfy the hypotheses of
Corollary~\ref{cor:uniform} and attain the lower bound \(n\).
The index formula also gives \(i(y_j,m)\equiv n\pmod2\), so
all iterates are good and \(\mu_{y_j}(m)\equiv n+1\pmod2\).
Their variational degrees are therefore even, and
\(\mchi(y_j)=1\). The example establishes sharpness of the
multiplicity bound. Since every \(\Delta_{y_j}>2\), it does not
establish optimality of the mean index condition.
\end{example}

\section{A smaller bound with an index exclusion}\label{sec:smaller}

Throughout this section, \(n\ge3\) is odd. We now prove
Theorem~\ref{thm:smaller}. The exclusion of good
iterates of Conley--Zehnder index zero gives an additional restriction
on the orbit counts. To obtain that restriction, we use two common
index jumps whose offsets have opposite signs. The coefficient bound
and the asymptotic counting formula remain those proved above.

Retain the notation and the assumptions of
Lemma~\ref{lem:recurrence}. Two sets of jump data
\[
 J=\bigl(d,(k_y,s_y)_{y\in\PP}\bigr),\qquad
 J'=\bigl(d',(k'_y,s'_y)_{y\in\PP}\bigr)
\]
of the same order are called \emph{complementary} if
\(s'_y=-s_y\) for every \(y\in\PP\). We call \(J\)
\emph{normalized} if
\(\sum_{y\in\PP}k_y\mchi(y)=d/2\).
The complementary jump data used below can also be obtained from
the generalized common index jump theorem and the construction using opposite
vertices described in
\cite[Theorem~3.6, Remark~3.7, and Step~2 in Case~1 of the proof
of Theorem~1.2]{DLLW}. The normalization follows from the
resonance identity \eqref{eq:resonance} and the integrality
argument used in Claim~1 of the same proof. We give a direct
construction in our notation.

\begin{lemma}\label{lem:paired}
Let \(n\ge3\) be odd. Assume that every iterate is nondegenerate,
that \(\PP\) is finite, that \(\Delta_y>0\) for every
\(y\in\PP\), and, in addition, that \(\mu_y(m)\ne0\) whenever
\(m\in G_y\). For every \(L_0,H\in\N\), there exist complementary
normalized jump data \(J,J'\) of order at least \(L_0\), with
\(d,d'\ge H\), such that
\begin{equation}\label{eq:paired}
 a+z(J)-c\ge n+1,\qquad 2a-c\ge n+1,
 \qquad z(J)=\#\{y\in\PP_1:s_y=0\}.
\end{equation}
\end{lemma}

\begin{proof}
Choose an integer \(L\ge L_0\) such that
\(L>3(n-1)/\delta_*\), and put \(H_* =\max\{H,n+3\}\).
We first construct the two normalized sets of jump data. We then
compare the alternating Morse sums below their respective centers.

\smallskip
\noindent\emph{Construction of the two jumps.}
Use the construction in the proof of Lemma~\ref{lem:recurrence}
with order \(L\), center at least \(H_*\), and approximation error
\(\varepsilon>0\) satisfying the conditions used there and also
\[
 2\varepsilon\sum_{y\in\PP}|\mchi(y)|<1.
\]
In the notation of that proof, the resulting
integers \(N,k_y,b_{yj}\), with \(d=2N\), satisfy
\[
 u_y=\frac{N}{\Delta_y}-\frac{k_y}{2},\qquad |u_y|<\varepsilon,
 \qquad
 e_{yj}=k_y\alpha_{yj}-b_{yj},\qquad |e_{yj}|<3\varepsilon,
\]
as well as \eqref{eq:integer-center}. Moreover,
\[
 s_y=\sum_j w_{yj}\operatorname{sgn}(e_{yj}).
\]
Here each \(e_{yj}\) is nonzero because \(\alpha_{yj}\) is
irrational and \(k_y>0\).

Multiplying \eqref{eq:resonance} by \(2N\) and using
\(k_y=2N/\Delta_y-2u_y\), we obtain
\[
 \left|\sum_{y\in\PP}k_y\mchi(y)-N\right|
 \le 2\varepsilon\sum_{y\in\PP}|\mchi(y)|<1.
\]
The expression inside the absolute value is an integer: every
\(k_y\) is even, and \(2\mchi(y)\in\Z\). Consequently
\(\sum_y k_y\mchi(y)=N\).

Choose \(0<\varepsilon_0\le\varepsilon\) such that
\(3\varepsilon_0<|e_{yj}|\) for all angular terms. Apply the
same approximation construction once more, with error
\(\varepsilon_0\), and denote the resulting integers and errors by
\(N_0,k_{0,y},b_{0,yj}\) and
\(e_{0,yj}=k_{0,y}\alpha_{yj}-b_{0,yj}\).
Since the approximation gives arbitrarily large \(N_0\), we may
require
\[
 2(N_0-N)\ge H_*,\qquad k_{0,y}-k_y>L
 \quad(y\in\PP).
\]
The normalization just proved also applies to these data, since
\(\varepsilon_0\le\varepsilon\). Set
\[
 N'=N_0-N,\qquad d'=2N',\qquad
 k'_y=k_{0,y}-k_y,\qquad b'_{yj}=b_{0,yj}-b_{yj}.
\]
Subtracting the corresponding identities for the two approximations
gives
\begin{equation}\label{eq:paired-normalization}
 k'_yA_y+2\sum_j w_{yj}b'_{yj}=d',\qquad
 \sum_y k'_y\mchi(y)=N'.
\end{equation}
The angular error
\[
 e'_{yj}=k'_y\alpha_{yj}-b'_{yj}=e_{0,yj}-e_{yj}
\]
has sign opposite to \(e_{yj}\), because
\(|e_{0,yj}|<3\varepsilon_0<|e_{yj}|\). Furthermore,
\(|e'_{yj}|<2|e_{yj}|<6\varepsilon\).
Thus addition of \(e'_{yj}\) leaves the integer parts of
\(\ell\alpha_{yj}\) and \(-\ell\alpha_{yj}\) unchanged for
\(1\le\ell\le L\). The computation in
Lemma~\ref{lem:recurrence}, using
\eqref{eq:paired-normalization}, now gives
\[
 \mu_y(k'_y\pm\ell)=d'\pm\mu_y(\ell),\qquad
 \mu_y(k'_y)-d'
   =\sum_j w_{yj}\operatorname{sgn}(e'_{yj})=-s_y.
\]
For \(t_y=0\), both offsets are zero. Hence the resulting data
\(J,J'\) are complementary and normalized, and both centers are
at least \(H_*\). Estimate~\eqref{eq:far} also holds with
\(d',k'_y\) in place of \(d,k_y\).

\smallskip
\noindent\emph{Alternating sums below the centers.}
Define
\[
 P(J)=\#\{y\in\PP_1:s_y\ge0\}
      -\#\{y\in\PP_{-1}:s_y\ge0\},\qquad
 P_+=P(J),\quad P_-=P(J').
\]
The condition defining \(P_-\) is \(s'_y\ge0\), or equivalently
\(s_y\le0\). For \(r\in\Z\), let
\(S(r)=\sum_{j<r}(-1)^jC_j\); this sum is finite because the
coefficients vanish in all sufficiently negative degrees.
By \eqref{eq:complete-count} and the normalization of \(J\),
\[
 \sum_{y\in\PP}\ 
 \sum_{\substack{m\in G_y\\1\le m\le k_y}}
 (-1)^{\mu_y(m)}
 =\sum_{y\in\PP}k_y\mchi(y)=N.
\]
We compare this sum, truncated by iteration number, with \(S(d)\),
which is truncated by degree. A good central iterate \(y^{k_y}\)
is included in the first sum and excluded from \(S(d)\) exactly
when \(s_y\ge0\). Since \(k_y\) is even, its sign is positive
on \(\PP_1\) and negative on \(\PP_{-1}\). The total central
correction is therefore \(-P_+\).

If \(|m-k_y|>L\), estimate \eqref{eq:far} shows that
\(\mu_y(m)<d\) exactly when \(m<k_y\), so there is no
correction from these iterates. For \(1\le\ell\le L\), put
\(r=\mu_y(\ell)\). The two iterates of orders \(k_y-\ell\)
and \(k_y+\ell\) have degrees \(d-r\) and \(d+r\),
respectively. Both are good exactly when \(\ell\in G_y\),
and their signs in the alternating sum agree.
If \(r>0\), both truncations select the iterate of order
\(k_y-\ell\). If \(r<0\), replacing that iterate by the one
of order \(k_y+\ell\) does not change the signed count.
If \(r=0\), neither degree is less than \(d\), whereas the
iteration truncation includes \(k_y-\ell\); the correction is
then \(-1\). Every good iterate of degree zero lies in this
range, since \eqref{eq:error} gives
\(\ell\delta_*\le n-1<L\delta_*\).
It follows that
\begin{equation}\label{eq:paired-truncation}
 S(d)=N-P_+-C_0,\qquad
 S(d')=N'-P_--C_0.
\end{equation}
The second identity is the same comparison for \(J'\).

Because \(d=2N\ge n+3\), the nonzero Betti coefficients below
\(d\) occur in degrees \(n+1,n+3,\ldots,d-2\), and their
sum is \(N-(n+1)/2\). In the alternating sum of
\(U_j+U_{j-1}\) over \(j<d\), all terms cancel except
\(-U_{d-1}\). Thus the Morse relation \eqref{eq:morse} and
\eqref{eq:paired-truncation} give
\begin{equation}\label{eq:paired-exact}
 P_++C_0=\frac{n+1}{2}+U_{d-1},\qquad
 P_-+C_0=\frac{n+1}{2}+U_{d'-1}.
\end{equation}
The additional index exclusion gives \(C_0=0\). Hence
\(P_++P_-\ge n+1\).
An orbit in \(\PP_1\) enters this sum once if \(s_y\ne0\),
and twice if \(s_y=0\). An orbit in \(\PP_{-1}\) has odd
degree \(\mu_y(k_y)=d+s_y\), so \(s_y\ne0\), and is
subtracted exactly once. Therefore
\(P_++P_-=a+z(J)-c\). Since \(z(J)\le a\), both
inequalities in \eqref{eq:paired} follow.
\end{proof}

\begin{remark}
The complementary construction and the exact identities
\eqref{eq:paired-exact} require no index exclusion. The exclusion
is used only to set \(C_0=0\). In the variational grading of
Duan and Fan, \(M_{-n-1}=C_0\). Thus dropping the nonnegative
remainders in \eqref{eq:paired-exact} recovers
\cite[Lemma~4.1, (4.34)--(4.35)]{DF}, whose Section~4.1 treats
general odd \(n\). The count above records the zero offsets for
each such dimension.
\end{remark}

\begin{proof}[Proof of Theorem~\ref{thm:smaller}]
Suppose that \(q\le n-1\), so that \(\PP\) is finite.
By Lemma~\ref{lem:paired}, \(2a-c\ge n+1\), whereas
\(a+c\le q\le n-1\). Therefore
\[
 3(a-c)=2(2a-c)-(a+c)\ge n+3.
\]
With \(D=n-a+c\), it follows that
\begin{equation}\label{eq:smaller-range}
 1\le D\le D_{\max}(n)
 :=\left\lfloor\frac{2n-3}{3}\right\rfloor.
\end{equation}
Also,
\[
 2c+h=q-b-a+c\le n-1-a+c=D-1.
\]

By \eqref{eq:rho-lower} in Lemma~\ref{lem:counting},
we have \(2\rho\ge D\ge1\), and hence \(\rho>0\).
It follows from \eqref{eq:rho} that the set
\[
 \mathcal O
 :=\{y\in\PP:\mchi(y)<0\}
 =\PP_{-1}\cup\PP_{-1/2}
\]
is nonempty. Define
\[
 \delta=\min_{y\in\mathcal O}\Delta_y>0.
\]
For every \(y\in\mathcal O\), the mean index assumption in
Theorem~\ref{thm:smaller} gives \(\Delta_y>\theta_n\).
Since \(\mathcal O\) is finite and nonempty, we obtain
\[
 \delta>\theta_n=1-\frac1{D_{\max}(n)}.
\]
Moreover, \(\Delta_y\ge\delta\) for every
\(y\in\PP_{-1}\cup\PP_{-1/2}\). Thus
\eqref{eq:rho-upper} in Lemma~\ref{lem:counting} applies.
Combining this upper bound with \eqref{eq:rho-lower}
and the inequality \(2c+h\le D-1\) proved above yields
\begin{equation}\label{eq:smaller-squeeze}
 D\le2\rho\le\frac{2c+h}{\delta}  \le\frac{D-1}{\delta}.
\end{equation}

On the other hand, \eqref{eq:smaller-range} gives
\[
 \frac{D-1}{D}\le1-\frac1{D_{\max}(n)}<\delta,
\]
and hence \((D-1)/\delta<D\), contradicting
\eqref{eq:smaller-squeeze}. Therefore \(q\ge n\).
For \(n=3\), one has \(D_{\max}(3)=1\) and \(\theta_3=0\),
so the strict mean index condition is already contained in the
positivity assumption.
\end{proof}

\begin{remark}[The critical value]\label{rem:critical}
For \(n\ge5\), the preceding argument does not include the
critical value \(\delta=\theta_n\).
If only \(\Delta_y\ge\theta_n\) is assumed on
\(\PP_{-1}\cup\PP_{-1/2}\), a hypothetical configuration with \(q\le n-1\)
would force equality throughout \eqref{eq:smaller-squeeze},
with \(D=D_{\max}(n)\) and \(2\rho=D\).

To see why the coefficient argument does not then give a
contradiction, put
\(\nu_n(T)=\#([n,T]\cap(2\Z+1))=T/2+O(1)\).
The bounded remainder in \eqref{eq:asymptotic} gives
\[
 \sum_{\substack{n\le p\le T\\p\ {\rm odd}}}(C_p-D)
 =A(T)-D\nu_n(T)=O(1).
\]
Each summand is a nonnegative integer. Hence \(C_p=D\) for
all sufficiently large odd \(p\). Since
\(C_p=U_{p-1}+U_p\) and \(U_{p-1}\ge D\), one has
\(U_{p-1}=D\) and \(U_p=0\) for those \(p\). The Morse
relation consequently gives \(C_j=D+1\) in every sufficiently
large even degree.

Choose \(L>(n+3(n-1))/\delta_*\) and take the centers in
Lemma~\ref{lem:paired} sufficiently large that \(d\ge2n\)
and \(C_d=D+1\). Applying the reflection formula with \(p=n\)
gives \(C_d=2C_0+Z_0(J)=z(J)\). Here \(C_0=0\), and a good central
iterate of even degree belongs to \(\PP_1\).
Lemma~\ref{lem:paired} gives \(z(J)\ge n+1-a+c=D+1\),
which is consistent with \(C_d=D+1\).
Thus this comparison does not exclude the critical case.
We do not know whether the critical value is sufficient, nor
whether the sufficient bound in Theorem~\ref{thm:smaller} is optimal.
\end{remark}

\end{document}